\documentclass[12pt] {article}
\usepackage{amsmath,amsthm}
\usepackage{amssymb,latexsym}

\usepackage{empheq}
\usepackage{amscd}
\usepackage{chngcntr}
\usepackage{enumitem}
\newlist{paragraphlist}{enumerate}{1}

\setlist[paragraphlist,1]{leftmargin=*,label={\bfseries \arabic*}}

\counterwithin{paragraphlisti}{subsubsection}

\title{On convergence of intrinsic volumes of Riemannian manifolds.}
\date{}
\author{Semyon Alesker \footnote{Partially supported by ISF grant 865/16 and the  US - Israel BSF grant 2018115.}
\\  { \normalsize Department of Mathematics, Tel Aviv University, Ramat Aviv}
\\  { \normalsize 69978 Tel Aviv, Israel }
\\ {\normalsize e-mail: semyon@tauex.tau.ac.il }}

\def\RR{\mathbb{R}}

\def\eps{\varepsilon}
\def\alp{\alpha}
\def\ome{\omega}

\def\to{\longrightarrow}

\def\qed { Q.E.D. }

\def\inj{\hookrightarrow}

\swapnumbers
\newtheorem{theorem}{Theorem}[section]

\newtheorem{lemma}[theorem]{Lemma}
\newtheorem{proposition}[theorem]{Proposition}
\newtheorem{claim}[theorem]{Claim}
\theoremstyle{definition}

\newtheorem{definition}[theorem]{Definition}
\newtheorem{remark}[theorem]{Remark}
\theoremstyle{conjecture}

\theoremstyle{principle}

\def\pt{\partial}

\numberwithin{equation}{section}
\begin{document}
\maketitle

\begin{abstract}
Let $\pi \colon M\to B$ be a Riemannian submersion of two compact smooth Riemannian manifolds, $B$ is connected.
Let $M(\eps)$ denote the manifold $M$ equipped with the new Riemannian metric which is obtained from the original one by multiplying by $\eps$ along the vertical subspaces
(i.e. along the fibers) and keeping unchanged along the (orthogonal to them) horizontal subspaces. Let $V_i(M(\eps))$ denote the $i$th intrinsic volume. The main result of this note says that
$\lim_{\eps\to +0}V_i(M(\eps))=\chi(Z) V_i(B)$ where $\chi(Z)$ denotes the Euler characteristic of a fiber of $\pi$.
\end{abstract}

\section{Introduction.}\label{S:introduction}
\begin{paragraphlist}
\item Intrinsic volumes of convex subsets of a Euclidean space have been studied in convexity for a long time, see e.g. \S 4.1 of Schneider's book \cite{schneider-book}. In dimension 2 and 3 this notion goes back to Steiner (1819). 
They can be defined using the Steiner formula as follows.
Let $\RR^n$ denote the Euclidean space. Let $B\subset \RR^n$ denote the unit Euclidean ball. Let $K\subset \RR^n$ be an arbitrary convex compact subset. Set $K_\eps$ denotes $\eps$-neighborhood of $K$.
The Steiner formula says that volume of $K_\eps$ is a polynomial of degree $n$ in $\eps\geq 0$:
\begin{eqnarray}\label{E:Steiner}
vol(K_\eps)=\sum_{i=0}^n\eps^{n-i}\kappa_{n-i}V_i(K),
\end{eqnarray}
where $\kappa_{n-i}$ are normalizing constants (see \cite{schneider-book} for an explicit expression). For $i=0,1,\dots, n$ the coefficient $V_i(K)$ is called the $i$th intrinsic volume of $K$. Note that $V_n(K)$ equals the volume of $K$ while $V_0(K)=1$.
Among important properties of intrinsic volumes of convex bodies let us mention their invariance under isometries of $\RR^n$, continuity in the Hausdorff metric, and the valuation property (i.e. they are finitely additive measures on convex compact sets); see \cite{schneider-book} for details. 
Intrinsic volumes are distinguished by the famous Hadwiger's characterization theorem \cite{hadwiger-thm} saying that a functional on all convex compact subsets of $\RR^n$ is a linear combination of intrinsic volumes if and only if it is
invariant under Euclidean isometries, continuous in the Hausdorff metric, and is a valuation.  Besides that, intrinsic volumes satisfy the Alexandrov-Fenchel inequality (see Ch. 7 in \cite{schneider-book}).

\item In 1939 H. Weyl \cite{weyl-39} has extended the definition of intrinsic volumes to smooth closed Riemannian manifolds. In this context intrinsic volumes are often called Lipschitz-Killing curvatures.
This notion is central for the present paper. The definition is based on a non-trivial result of him stated below.
However we will state it in greater generality
of manifolds with boundary or, still more generally, simplicial corners. Such an extension seems to be a folklore, and this generality seems to be more natural to us.

 Let $(M^n,g)$ be a compact smooth Riemannian manifold, possibly with simplicial corners (e.g. with boundary).  By the Nash imbedding theorem there exists a smooth isometric imbedding $\iota\colon M\inj \RR^L$ into 
 a Euclidean space of sufficiently high dimension. The Euclidean volume  of the $\eps$-neighborhood of the image of $\iota$ turns out to be polynomial for small $\eps\geq 0$:
 $$  vol_L(\iota M)_\eps=\eps^{L-n}\sum_{i=0}^n c_{L,n,i} \eps^{n-i} V_i(M),$$
where $c_{L,n,i}$ are appropriately normalized coefficients independent of $M$ and $\iota$. The Weyl theorem says that the coefficients $V_i(M)$ are independent of the imbedding $\iota$. Weyl himself proved it for closed manifolds; he gave an explicit formula 
(see below) in terms of the Riemann curvature tensor.
The case of boundary and corners seems to be known for a long time, but I do not have a classical reference. But this is a special case of Theorem 3.11 in a recent paper \cite{fu-wannerer}. This definition is obviously consistent with the 
definition for convex sets from the previous paragraph when convex sets have smooth boundary.

Note that $V_n(M^n)$ equals to the Riemannian volume of $M$, $V_0(M)$ equals to the Euler characteristic of $M$ by the generalized Gauss-Bonnet formula discussed below. 
For a closed $M$ one has $V_i(M^n)=0$ provided $n-i$ is odd. Moreover for  closed $M$ the intrinsic volume $V_{n-2}(M)$ is proportional
to the integral of the scalar curvature $\int_M Sc(x) dvol_M(x)$.

\item As we have mentioned, an important property of intrinsic volumes on convex sets is continuity in the Hausdorff metric. The convergence properties of intrinsic volumes of manifolds seem to be rare in literature.
The earliest known to me result in this direction is due to Cheeger, M\"uller, and Schrader \cite{ch-mu-sch} who proved a convergence of intrinsic volumes for polyhedral approximations of smooth Riemannian manifolds.

A more recent example is the work in progress by Lebedeva and Petrunin \cite{lebedeva-petrunin} where it is claimed that for a non-collapsing sequence $\{M_i^n\}$ of closed Riemannian manifolds
with uniformly bounded below sectional curvature converging in the Gromov-Hausdorff sense to an Alexandrov space always there exists $\lim_{i\to \infty} V_{n-2}(M^n_i)$.

The author (based on discussions with A. Petrunin \cite{petrunin-personal}) has formulated a few conjectures \cite{alesker-conjectures} on the behavior of intrinsic volumes under the Gromov-Hausdorff convergence 
of Riemannian manifolds, and more generally Alexandrov spaces, with sectional curvature uniformly bounded below. The main result of the present paper verifies one of these conjectures in a very special 
case. In fact in this case the assumption of uniform lower boundedness of sectional curvature is unnecessary. 

\item Let us formulate the main result. 
\begin{definition}\label{rieman-submersion}
A smooth map $\pi\colon M\to B$ of two smooth Riemannian manifolds is called a Riemannian submersion if for any $o\in M$ the differential 
$$d\pi_o\colon T_oM\to T_{\pi(o)}B$$
being restricted to the orthogonal complement of its kernel $(Ker (d\pi_o))^\perp$ is an isometry.
\end{definition}
Note that $Ker (d\pi_o)$ is the tangent space $T_oZ_o$ of the fiber $Z_o:=\pi^{-1}\pi(o)$.

Let $\pi\colon M\to B$ be a Riemannian submersion. Let us change the original Riemannian metric $g$ on $M$ as follows. Fix any $o\in M$. We have the orthogonal with respect to $g$ decomposition of the tangent space
$$T_oM=T_oZ_o\oplus (T_oZ_o)^\perp.$$
For $\eps >0$ let us denote by $g(\eps)$ the metric which is equal to $\eps g$ on $T_oZ_o$, $g$ on $(T_oZ_o)^\perp$, and the two subspaces are still orthogonal to each other with respect to $g(\eps)$. We denote by $M(\eps)$
the manifold $M$ equipped with the metric $g(\eps)$. The main result is 
\begin{theorem}\label{T:main-result}
Let $\pi\colon M\to B$ be a Riemannian submersion of smooth compact manifolds. Let us assume that $B$ is connected. Let $Z$ denote a fiber of $\pi$. Then for any $i\geq 0$ 
$$\lim_{\eps\to +0}V_i(M(\eps))=\chi(Z) V_i(B),$$
where $\chi(Z)$ denotes the Euler characteristic of $Z$.
\end{theorem}
Just to compare with the conjectures \cite{alesker-conjectures}, we show in Proposition \ref{P:sect-curvature} that under the assumptions of Theorem \ref{T:main-result} the manifolds $M(\eps)$ have uniform lower bound on sectional curvature
when $\eps\to +0$ if and only if all the fibers of $\pi$ have non-negative sectional curvature.

\item Intrinsic volumes on convex sets and on manifolds play an important role in valuations theory on convex sets \cite{alesker-kent}, \cite{schneider-book} and its generalizations to manifolds \cite{alesker-mflds1}-\cite{alesker-mflds4}, \cite{fu-wannerer}.
They play a key role in integral geometry \cite{bernig-fu-unitary}, \cite{bernig-fu-solanes-unitary}, \cite{klain-rota}.

Recently there were constructed analogues of intrinsic volumes in pseudo-Riemannian geometry and obtained non-trivial results for them \cite{bernig-faifman}, \cite{bernig-faifman-solanes}, \cite{bernig-faifman-solanes2}.

\item {\bf Acknowledgement.} Part of this work was done during my sabbatical stay at the Kent State University in the academic year 2018/19. I am grateful to this institution for hospitality.
\end{paragraphlist}

\section{Reminder on Weyl's and generalized Gauss-Bonnet formulas.}\label{S:background}
First let us remind the H. Weyl \cite{weyl-39} formula for intrinsic volumes of a closed smooth Riemannian manifold $M^n$. 
Let $R_{abcd}$ denote the Riemann curvature tensor for a Riemannian metric $g_{ab}$.
Let us denote
$$R^{pq}_{rs}:=g^{pa}g^{qb}R_{abrs},$$
where $(g^{ab})$ denotes the inverse matrix of $(g_{ab})$, and we use the summation convention over repeated indices.

H. Weyl \cite{weyl-39} has shown that for a closed manifold $M^n$ for odd $e$ one has $V_{n-e}(M^n)=0$, and for even $e$ one has
\begin{eqnarray}\label{E:IntrinVol}
V_{n-e}(M^n)=(2\pi)^{-e/2}\int_M \sum_{[p,q]}sgn(p,q)R^{q_1q_2}_{p_1p_2}\dots R^{q_{e-1}q_e}_{p_{e-1}p_e}dvol_M,
\end{eqnarray}
where $dvol_M$ is the Riemannian volume form and the sum runs over all couplings $[p,q]$. Here the numbers $p_1,\dots,p_e$ are distinct numbers between 1 and $n$, and $q_1,\dots,q_e$ are any permutation of the $p_i$'s.
The sign $sgn(p,q)$ is the sign of the permutation taking $p_1,\dots,p_e$ to $q_1,\dots,q_e$.
The coupling $[p,q]$ is a partition of pairs
\begin{eqnarray*}
\begin{array}{c|c|c}
(p_1p_2)&(p_3p_4)&\dots\\
(q_1q_2)&(q_3q_4)&\dots
\end{array}.
\end{eqnarray*}
The pairs $(p_1p_2)$ etc. of two distinct numbers is unordered; due to symmetries of the curvature tensor and the term $sgn(p,q)$ the switch of the order of 
$p_1$ and $p_2$ does not change the summand. The coupling is also unchanged under any permutation of any of $e/2$ columns of it.

\hfill

In 1943 Allendoerfer and Weil \cite{allendoerfer-weil} have shown that $V_0(M)$ defined by the formula (\ref{E:IntrinVol}) (with $e=n$, $n$ is even) coincides with the Euler characteristic $\chi(M)$ thus obtaining a generalization of the Gauss-Bonnet formula. 
Actually they got a more general result for compact manifolds with simplicial corners. In 1944 Chern \cite{chern-1944} published a simpler proof of this fact for the case of closed manifolds.

\section{A few technical lemmas.}\label{S:technical-lemmas}
Let $\pi\colon M\to B$ be a Riemannian submersion of compact manifolds. We fix a point $o\in M$. Define $o^*:=\pi(o)$. There exists a neighborhood $U$ of $o^*$ such that $\pi^{-1}(U)$ is diffeomorphic to $Z\times U$ and 
under this identification $\pi$ is the obvious projection $Z\times U\to U$ where $Z$ is a fiber of $\pi$. In all local computations throughout this paper we will fix an arbitrary such a diffeomorphism.

We fix coordinates $\{x^\alp\}$ near $o^*$.
Let us denote $Z_{o^*}:=\pi^{-1}(o^*)$. Let us fix coordinates $x^i$ on $Z_{o^*}$ near $o$. Throughout this paper we denote by Greek letters $\alp,\beta,\gamma,\delta\dots$ the coordinates on an open subset $U\subset B$, by Latin letters $i,j,k,l, m,n$
coordinates on the fiber $Z_{o^*}$, and by Latin letters $p,q,r,s,t$ coordinates of either kind on $Z\times U$.
Let $(g(\eps))$ the metric on $M$ as in Theorem \ref{T:main-result}. 
\begin{lemma}\label{L:new-metric}
At any point $o\in M$ we have
\begin{eqnarray*}
(g(\eps))=\left[\begin{array}{cc}
                \eps\cdot (g_{ij})&\eps\cdot (g_{i\alp})\\
                \eps\cdot (g_{\alp i})&g_{\alp\beta}
                \end{array}\right].
                \end{eqnarray*}
\end{lemma}
{\bf Proof.} At any point $o\in M$ we denote by $(\cdot,\cdot)$ the scalar product on $T_0M$ for the original metric $g$, and by $(\cdot,\cdot)_\eps$ for the metric $g(\eps)$.

First it is obvious that $g(\eps)_{ij}=\eps\cdot g_{ij}$. 
Next in the orthogonal complement to the tangent space at $o$ to the fiber there exists a unique basis $\{\xi_\alp\}$
such that $d\pi(\xi_\alp)=\pt_\alp$ where $\pt_\alp$ are the obvious vector fields corresponding to chosen coordinate $x^\alp$.
Clearly 
\begin{eqnarray}\label{E:scalar-product}
g(\eps)_{\alp\beta}=(\pt_\alp,\pt_\beta)_\eps=(\xi_\alp,\xi_\beta)_\eps=(\xi_\alp,\xi_\beta)=g_{\alp\beta},
\end{eqnarray}
where the second equality follows from the Riemannian submersion property of $\pi$.

It remains to prove $g(\eps)_{i\alp}=\eps g_{i\alp}$. Let us denote by $ (\tilde g^{ij})$ the matrix inverse to $(g_{ij})$ (note that both are matrices of size $N\times N$ where $N$ is the dimension of the fiber of $\pi$).
First let us express $\xi_\alp$ via $\pt_p$'s. Since $d\pi(\xi_\alp)=\pt_\alp$ we have 
$$\xi_\alp=\pt_{\alp}+h^i_\alp \pt_i,$$
where $h^i_\alp$ are some coefficients. Then 
\begin{eqnarray*}
0=(\xi_\alp,\pt_j)=(\pt_{\alp},\pt_j)+h^i_\alp(\pt_i,\pt_j)=g_{\alp j}+h^i_\alp g_{ij}.
\end{eqnarray*}
Consequently $h^i_\alp=-g_{\alp j}\tilde g^{ji}.$
Hence 
\begin{eqnarray}\label{E:xi-pt}
\xi_\alp=\pt_{\alp}-g_{\alp j}\tilde g^{ji}\pt_i.
\end{eqnarray}

Using this we have
 
\begin{eqnarray*}
g(\eps)_{i\alp}=(\pt_i,\pt_\alp)_\eps=(\pt_i,\xi_\alp+g_{\alp j}\tilde g^{jk}\pt_k)_\eps=g_{\alp j}\tilde g^{jk}(\pt_i,\pt_k)_\eps=\\
\eps g_{\alp j}\tilde g^{jk}g_{ik}=\eps g_{\alp i}=\eps g_{i\alp}.
\end{eqnarray*}
\qed

\begin{lemma}\label{L:limit-space}
In the above notation the Riemannian metric on $B$ at $o^*=\pi(o)$ is equal to $(g_{\alp\beta})$.
\end{lemma}
{\bf Proof} immediately follows from \ref{E:scalar-product}. \qed

\begin{lemma}\label{L:estimate-Christoffel}
Let us choose the normal coordinate system near $o^*=\pi(o)$, (in particular $\pt_\gamma g_{\alp\beta}|_o=0$).
Then as $\eps\to +0$ one has
\begin{eqnarray*}
\pt_p g_{qr}(\eps)|_o=O(\eps),\\
\Gamma_{pqr}(\eps)|_o=O(\eps).
\end{eqnarray*}
\end{lemma}
{\bf Proof.} Obviously the second part follows from the first one since 
$$\Gamma_{pqr}(\eps)=\frac{1}{2}(\pt_r g_{pq}(\eps) +\pt_q g_{pr}(\eps)-\pt_p g_{qr}(\eps)).$$ 
To prove the first part observe that if at least one of $q,r$ is $i,j,...$ then
$g_{qr}(\eps)=\eps g_{qr}$, and the result follows in this case. In the remaining case $q=\alp,r=\beta$ we have
$$\pt_p g_{\alp \beta}(\eps)|_o=\pt_pg_{\alp\beta}|_o.$$
By Lemma \ref{L:limit-space} $g_{\alp\beta}$ is constant along the fibers of $\pi$. Hence $\pt_{i}g_{\alp\beta}|_o=0$.
By the choice of normal coordinate system near $o^*$ we have $\pt_{\gamma}g_{\alp\beta}|_o=0$. \qed

\begin{lemma}\label{L:matrices}
Let \begin{eqnarray*}
    G=\left[\begin{array}{cc}
              A& B\\
              C&D
            \end{array}\right]
            \end{eqnarray*}
           be a block matrix with $A,D$ invertible square matrices of not necessarily the same size. Then for $\eps \to 0$ one has
\begin{eqnarray*}
\left[\begin{array}{cc}
              \eps A&\eps B\\
             \eps C& D
            \end{array}\right]^{-1}=\left[\begin{array}{cc}
                                     \eps^{-1} A^{-1}+O(1)&O(1)\\
                                     O(1)&D^{-1}+O(\eps)
                                    \end{array}\right]
\end{eqnarray*}
while the left hand side is automatically invertible for small $\eps>0$.
\end{lemma}
{\bf Proof.} We have
\begin{eqnarray*}
\left[\begin{array}{cc}
              \eps A&\eps B\\
             \eps C& D
            \end{array}\right]=\left[\begin{array}{cc}
                                         A&0\\
                                         0&D
                                         \end{array}\right]
            \left[\begin{array}{cc}
                                      \eps I&\eps X\\
                                      \eps Y&I
                                     \end{array}\right],
\end{eqnarray*}
where $X=A^{-1}B$, $Y=D^{-1}C$. Hence it suffices to prove the lemma for $A=I,D=I$. It is easy to see that
\begin{eqnarray*}
\left[\begin{array}{cc}
              \eps I&\eps X\\
             \eps Y& I
            \end{array}\right]=\left[\begin{array}{cc}
                                   I&\eps X\\
                                   0&I
                                      \end{array}\right]\left[\begin{array}{cc}
                                                              \eps(I-\eps XY)&0\\
                                                              0&I
                                                               \end{array}\right]\left[\begin{array}{cc}
                                                                                          I&0\\
                                                                                          \eps Y&I
                                                                                         \end{array}\right].
\end{eqnarray*}
Taking the inverse matrix and computing the inverse of each term we finally obtain
\begin{eqnarray*}
\left[\begin{array}{cc}
              \eps I&\eps X\\
             \eps Y& I
            \end{array}\right]^{-1}=\left[\begin{array}{cc}
                                        \eps^{-1}(I-\eps XY)^{-1}&-(I-\eps XY)^{-1}X\\
                                        -Y(I-\eps XY)^{-1}&I+\eps Y(I-\eps XY)^{-1}X
                                      \end{array}\right].
\end{eqnarray*}
The lemma follows. \qed

\begin{lemma}\label{L:curvature-asymp}
At any point $o\in M$ one has
\begin{eqnarray}
R_{\alp\beta\gamma\delta}(\eps)=O(1), \label{E:curv2}\\\label{E:curv1}
R_{ipqr}(\eps)=O(\eps).
\end{eqnarray}
where the estimates in $O(\eps), O(1)$ depend on $||g_{pq}||_{C^2}$ only.
\end{lemma}
\begin{remark}
The lemma says that if at least one of the indices is not Greek then $R(\eps)=(\eps)$, otherwise it is just bounded. Indeed this follows from the symmetries of the Riemann curvature tensor.
\end{remark}

{\bf Proof.} Let us choose near $o^*=\pi(o)$ normal coordinates. There is a known formula for the Riemann curvature tensor (see e.g. formula (92.1) in \cite{landau-lifschitz}):
\begin{eqnarray}\label{E:RC-formula}
2R_{pqrs}(\eps)=\\\label{E:RC-formula2}
\pt_q\pt_r g_{ps}(\eps)+\pt_p\pt_s g_{qr}(\eps)-\pt_q\pt_s g_{pr}(\eps)-\pt_p\pt_r g_{qs}(\eps)+2(\Gamma_{tqr}(\eps)\Gamma^t_{ps}(\eps)-\Gamma_{tqs}(\eps)\Gamma^t_{pr}(\eps)).
\end{eqnarray}

Let us show that each of the two $\Gamma\Gamma$ summands is $O(\eps)$. Let us estimate the first one, the second one is estimated similarly.
We have
\begin{eqnarray*}
\Gamma_{tqr}(\eps)\Gamma^t_{ps}(\eps)=\Gamma_{tqr}(\eps)g^{tt'}(\eps) \Gamma_{t'ps}(\eps)=O(\eps),
\end{eqnarray*}
where we have used that 
by Lemma \ref{L:estimate-Christoffel} $\Gamma(\eps)=O(\eps)$ and by Lemmas \ref{L:new-metric} and \ref{L:matrices} $g^{tt'}(\eps)=O(\eps^{-1})$.

%

 If at least one of $p,q$ is $i,j,k,l$ then $g_{pq}(\eps)=O(\eps)$ and all of its derivatives. Hence in (\ref{E:RC-formula2}) the only contributions which may not be $O(\eps)$
 have the form $\pt_p\pt_q g_{\alp\beta}(\eps)$. By Lemma \ref{L:limit-space} (with $(g)$ replaced with $(g(\eps))$) we have
 $$\pt_i(g_{\alp\beta}(\eps))=\pt_i g_{\alp\beta}=0.$$
 Differentiating by $\pt_p$ (recall that $p$ corresponds to any coordinate on $M$) we get
 $$\pt_p\pt_i(g_{\alp\beta}(\eps))=0.$$
 This implies the lemma. \qed

 \hfill


\section{Proof of Theorem \ref{T:main-result}.}\label{S:proof1}
Let us apply (\ref{E:IntrinVol}) to $(M,g(\eps))=:M(\eps)$. We denote by $N$ the dimension of the fibers of $\pi$. We have
\begin{eqnarray}\label{E:intrVol-eps}
V_{n-e}(M(\eps))=(2\pi)^{-e/2}\eps^{N/2}\int_M\sum_{[p,q]}sgn(p,q)R^{q_1q_2}_{p_1p_2}(\eps)\dots R^{q_{e-1}q_e}_{p_{e-1}p_e}(\eps)dvol_M,
\end{eqnarray}
where $dvol_M$ is the Riemannian volume form corresponding to the original metric $g$ (rather than $g(\eps)$).
Indeed it is easy to see that for the metric $g(\eps)$ the volume form is equal to $dvol_{M(\eps)}=\eps^{N/2}dvol_M$.

\begin{theorem}\label{T:main}
For a Riemannian submersion $\pi\colon M^n\to B$ of closed Riemannian manifolds
$$\lim_{\eps\to +0}V_i(M(\eps))=\chi(Z)\cdot V_i(B),$$
where $Z$ is a fiber of $\pi$.
\end{theorem}
{\bf Proof.} If $e$ is odd then $V_{n-e}(M(\eps))=0$ for all $\eps>0$. If $\dim B-n+e$ is odd then $V_{n-e}(B)=0$, and the theorem follows when both $e$ and $\dim B-n+e$ are odd.
If $e$ is odd and $\dim B-n+e$ is even then the dimension $N=n-\dim B$ of a fiber $Z$ is odd and hence $\chi(Z)=0$. Thus the result follows in this case too.

Thus let us assume that $e$ is even and consider $V_{n-e}(M(\eps))$ defined by (\ref{E:intrVol-eps}). By Lemmas \ref{L:curvature-asymp} and \ref{L:matrices} we have uniformly on $M$
\begin{eqnarray*}
R^{\alp\beta}_{\gamma\delta}(\eps)=O(1),\\
R^{\alp\beta}_{ip}(\eps)=O(\eps),\\
R^{i\alp}_{pq}(\eps)=O(1),\\
R^{ij}_{pq}(\eps)=O(\eps^{-1}).
\end{eqnarray*}
Let us consider a summand $R^{q_1q_2}_{p_1p_2}(\eps)\dots R^{q_{e-1}q_e}_{p_{e-1}p_e}(\eps)$ in (\ref{E:intrVol-eps}). Let us denote by
$\tau$ the number of terms of the form $R^{ij}_{pq}(\eps)$, and by $\ome$ the number of terms of the form $R^{\alp\beta}_{ip}(\eps)$. Then
the above formulas imply that uniformly on $M$
$$ R^{q_1q_2}_{p_1p_2}(\eps)\dots R^{q_{e-1}q_e}_{p_{e-1}p_e}(\eps)=O(\eps^{-\tau+\ome}).$$
Since all $p_1,\dots,p_e$ are pairwise distinct we have $\tau\leq N/2$. Hence uniformly on $M$
$$\eps^{N/2}R^{q_1q_2}_{p_1p_2}(\eps)\dots R^{q_{e-1}q_e}_{p_{e-1}p_e}(\eps)=O(\eps^{N/2-\tau+\ome}).$$
Thus the following claim holds.
\begin{claim}\label{Cl:0001}
The expression $\eps^{N/2}R^{q_1q_2}_{p_1p_2}(\eps)\dots R^{q_{e-1}q_e}_{p_{e-1}p_e}(\eps)=O(\eps^{N/2-\tau+\ome})$
\newline
(1) is always uniformly bounded on $M$ as $\eps\to +0$;
\newline
(2) vanishes in the limit when either $\tau<N/2$ or $\ome>0$.
\end{claim}
Hence it suffices to consider the summands with $\tau=N/2,\ome=0$ simultaneously. In particular $N$ must be even for such summands to exist.
Thus if $N$ is odd then $\lim_{\eps\to +0} V_{n-e}(M(\eps))=0$. At the same time $\chi(Z)=0$ and the theorem follows for odd $N$.

Let us assume that $N$ is even. Consider all the summands with $\tau=N/2$ and $\ome=0$.
For any of them all the  indices corresponding to coordinates of the fiber $Z$ appear among upper indices of the terms of the form $R^{ij}_{pq}(\eps)$, and
there are no more such indices to appear in terms of the form $R^{i\alp}_{pq}(\eps)$. In other words if $\tau=N/2$ there are no terms of the form $R^{i\alp}_{pq}(\eps)$.
Hence the only terms which appear are $R^{\alp\beta}_{\gamma\delta}(\eps)$ and $R^{ij}_{pq}(\eps)$. Furthermore since all lower indices are a permutation of upper indices,
the total numbers of Greek upper and lower indices should be the same. Hence the term $R^{ij}_{pq}(\eps)$ should be $R^{ij}_{kl}(\eps)$. Hence
\begin{eqnarray*}
(\star):=\lim_{\eps\to +0} \eps^{N/2} \sum_{[p,q]}sgn(p,q)R^{q_1q_2}_{p_1p_2}(\eps)\dots R^{q_{e-1}q_e}_{p_{e-1}p_e}(\eps)=\\
\lim_{\eps\to +0} \eps^{N/2}
\left(\sum_{[\alp,\beta]}sgn(\alp,\beta)R^{\alp_1\alp_2}_{\beta_1\beta_2}(\eps)\dots R^{\alp_{e-N-1}\alp_{e-N}}_{\beta_{e-N-1}\beta_{e-N}}(\eps)\right)\cdot\\
\left( \sum_{[i,j]}sgn(i,j)R^{i_1i_2}_{j_1j_2}(\eps)\dots R^{i_{N-1}i_{N}}_{j_{N-1}j_{N}}(\eps)\right).
\end{eqnarray*}
We claim that this limit exists and
\begin{eqnarray}\label{E:limit-important}
(\star)=\left(\sum_{[\alp,\beta]}sgn(\alp,\beta) \hat R^{\alp_1\alp_2}_{\beta_1\beta_2}\dots  \hat R^{\alp_{e-N-1}\alp_{e-N}}_{\beta_{e-N-1}\beta_{e-N}}\right)\cdot
\left( \sum_{[i,j]}sgn(i,j)\tilde R^{i_1i_2}_{j_1j_2}\dots\tilde R^{i_{N-1}i_{N}}_{j_{N-1}j_{N}}\right)
\end{eqnarray}
where $\hat R^{\alp\beta}_{\gamma\delta}$ denotes the curvature tensor of $B$ at the point $o^*$, and $\tilde R^{mn}_{kl}$ denotes the curvature tensor of the fiber $Z=\pi^{-1}\pi(o)$ at $o$.
To prove that it suffices to prove that following two statements:
\begin{eqnarray}\label{E:statement1}
\lim_{\eps\to +0}  R^{\alp\beta}_{\gamma\delta}(\eps) =\hat  R^{\alp\beta}_{\gamma\delta},\\\label{statement2}
\lim_{\eps\to +0}\eps R^{mn}_{kl}(\eps)=\tilde R^{mn}_{kl}.
\end{eqnarray}
For let us chose normal coordinates near $o^*$. First let us prove (\ref{E:statement1}). By (\ref{E:RC-formula}) one has at the point $o$
\begin{eqnarray}\label{E:line1}
2R^{\alp\beta}_{\gamma\delta}(\eps)=g^{\alp p}(\eps)g^{\beta q}(\eps)(\pt_\gamma\pt_q g_{\delta p}(\eps)+\pt_\delta\pt_p g_{\gamma q}(\eps)-\pt_{\gamma}\pt_p g_{\delta q}(\eps)-\pt_\delta\pt_q g_{\gamma p}(\eps)+\\\label{E:line2}
2(\Gamma_{t\delta p}(\eps)g^{ts}(\eps)\Gamma_{s\gamma q}(\eps)-\Gamma_{t\delta q}(\eps)g^{ts}(\eps)\Gamma_{s\gamma p}(\eps))).
\end{eqnarray}
By Lemma \ref{L:estimate-Christoffel} $\Gamma(\eps)=O(\eps)$. By Lemmas \ref{L:new-metric} and \ref{L:matrices} $g^{st}(\eps)=O(\eps^{-1})$. Hence the expression (\ref{E:line2}) tends to 0 when $\eps\to +0$.
By Lemma \ref{L:new-metric} $g_{\alp i}(\eps)$ and all of its derivatives are proportional to $\eps$. By Lemma \ref{L:matrices} $g^{\alp p}(\eps)=O(1)$. Consequently in (\ref{E:line1}) the summands with indices $p,q$ corresponding to coordinates of the fiber $Z=\pi^{-1}\pi o$
tend to 0. Hence
\begin{eqnarray*}
2R^{\alp\beta}_{\gamma\delta}(\eps)=g^{\alp \mu}(\eps)g^{\beta \nu}(\eps)(\pt_\gamma\pt_\nu g_{\delta \mu}(\eps)+\pt_\delta\pt_\mu g_{\gamma \nu}(\eps)-\pt_{\gamma}\pt_\mu g_{\delta \nu}(\eps)-\pt_\delta\pt_\nu g_{\gamma \mu}(\eps))+O(\eps)
\end{eqnarray*}
Let us denote by $\hat g_{\alp\beta}$ the Riemannian metric of the base $B$ at the point $o^*$. By Lemma \ref{L:new-metric} $g_{\alp\beta}(\eps)|_x=\hat g_{\alp\beta}|_{\pi(x)}$ for any $x$ in a neighborhood of $o$ and for any small $\eps$.
Hence $\pt_\gamma\pt_\delta g_{\alp\beta}(\eps)|_o=\pt_\gamma\pt_\delta \hat g_{\alp\beta}|_{o^*}$. By Lemma \ref{L:matrices}
$g^{\alp\beta}(\eps)=\hat g^{\alp\beta}+O(\eps)$. Hence
at the point $o$ one has
\begin{eqnarray*}
2R^{\alp\beta}_{\gamma\delta}(\eps)=\hat g^{\alp \mu}\hat g^{\beta \nu}(\pt_\gamma\pt_\nu \hat g_{\delta \mu}+\pt_\delta\pt_\mu \hat g_{\gamma \nu}-\pt_{\gamma}\pt_\mu \hat g_{\delta \nu}-\pt_\delta\pt_\nu \hat g_{\gamma \mu})+O(\eps)
\end{eqnarray*}
Since at $o^*$ we have chosen normal coordinates the Christoffel symbols of $\hat g$ vanish at $o^*$. By (\ref{E:RC-formula}) the right hand side of the last equality can be rewritten
$$2R^{\alp\beta}_{\gamma\delta}(\eps)=2\hat R^{\alp\beta}_{\gamma\delta} + O(\eps).$$
Thus (\ref{E:statement1}) follows. Let us now prove (\ref{statement2}). By (\ref{E:RC-formula}) we have (below as usual $m,n,k,l$ denote indices of coordinates along the fiber $Z$, and $p,q,r,s,t$ denote arbitrary indices)
\begin{eqnarray}\label{E:new-case1}
\eps 2R^{mn}_{kl}(\eps)=\eps g^{mp}(\eps)g^{nq}(\eps)(\pt_l\pt_p g_{kq}(\eps)+\pt_k\pt_qg_{lp}(\eps)-\pt_k\pt_pg_{lq}(\eps)-\pt_l\pt_qg_{kp}(\eps)+\\\label{E:new-case2}
2\left[\Gamma_{tlp}(\eps)g^{ts}(\eps)\Gamma_{skq}(\eps)-\Gamma_{tkp}(\eps)g^{ts}(\eps)\Gamma_{slq}(\eps)\right]).
\end{eqnarray}

By Lemma \ref{L:new-metric} each of the first summands in the round brackets is proportional to $\eps$, e.g.
\begin{eqnarray}\label{E:round1}
g_{kq}(\eps)=\eps g_{kq}.
\end{eqnarray}
By Lemmas \ref{L:new-metric}, \ref{L:matrices}, and \ref{L:estimate-Christoffel}
\begin{eqnarray*}
g^{mp}(\eps)=O(\eps^{-1}),\,\,\, \Gamma(\eps)=O(\eps).
\end{eqnarray*}
It follows that one opens brackets in (\ref{E:new-case1})-(\ref{E:new-case2}) then each summand is $O(1)$. Moreover if one of the terms in a summand has better estimate (e.g. $g^{m\alp}(\eps)=O(1)$ by Lemma \ref{L:matrices})
then the corresponding term is $O(\eps)$. Hence we may keep in (\ref{E:new-case1}) only terms of the form $g^{mi}(\eps)$ (i.e. both indices $m,i$ correspond to coordinates on $Z$). For the same reason $g^{st}(\eps)$
in (\ref{E:new-case2}) 
can be replaced with $g^{i'j'}(\eps)$ (i.e. summation is only over the indices corresponding to the fiber $Z$).
 Again by Lemma \ref{L:matrices}
$$g^{mi}(\eps)=\eps^{-1}(g|_Z)^{mi} +O(1).$$
Substituting all that into (\ref{E:new-case1})-(\ref{E:new-case2}) we get
\begin{eqnarray*}
\eps 2R^{mn}_{kl}(\eps)=\eps^{-1} (g|_Z)^{mi}(g|_Z)^{nj}(\eps(\pt_l\pt_i g_{kj}+\pt_k\pt_jg_{li}-\pt_k\pt_ig_{lj}-\pt_l\pt_jg_{ki})+\\
2\eps^{-1}\left[\Gamma_{i'li}(\eps)(g|_Z)^{i'j'}\Gamma_{j'kj}(\eps)-\Gamma_{i'ki}(\eps)(g|_Z)^{i'j'}\Gamma_{j'lj}(\eps)\right])+O(\eps).
\end{eqnarray*}
In the last formula all the indices correspond to coordinates of the fiber $Z$. For them obviously $\Gamma_{i'ki}(\eps)=\eps \Gamma_{i'ki}$ where $\Gamma_{i'ki}$ is the Christoffel symbol of the original metric on $M$.
Hence
\begin{eqnarray*}
\eps 2R^{mn}_{kl}(\eps)= (g|_Z)^{mi}(g|_Z)^{nj}(\pt_l\pt_i g_{kj}+\pt_k\pt_jg_{li}-\pt_k\pt_ig_{lj}-\pt_l\pt_jg_{ki})+\\
2\left[\Gamma_{i'li}(g|_Z)^{i'j'}\Gamma_{j'kj}-\Gamma_{i'ki}(g|_Z)^{i'j'}\Gamma_{j'lj}\right])+O(\eps)=\\
2\tilde R^{mn}_{ij}+O(\eps).
\end{eqnarray*}
Clearly statement (\ref{statement2}) follows. Consequently (\ref{E:limit-important}) is proven as well.

\hfill

Let us prove the theorem now.
By (\ref{E:intrVol-eps})
$$V_{n-e}(M(\eps))=(2\pi)^{-e/2}\int_M \eps^{N/2}\sum_{[p,q]}sgn(p,q)R^{q_1q_2}_{p_1p_2}(\eps)\dots R^{q_{e-1}q_e}_{p_{e-1}p_e}(\eps)dvol_M$$
By Claim \ref{Cl:0001}(1) the expression under the integral is uniformly bounded on $M$. By (\ref{E:limit-important}) and the Lebesgue convergence theorem
\begin{eqnarray*}
\lim_{\eps\to+0} V_{n-e}(M(\eps))=\\
(2\pi)^{-e/2}\int_M\left(\sum_{[\alp,\beta]}sgn(\alp,\beta) \hat R^{\alp_1\alp_2}_{\beta_1\beta_2}\dots  \hat R^{\alp_{e-N-1}\alp_{e-N}}_{\beta_{e-N-1}\beta_{e-N}}\right)|_{o}\cdot\\
\left( \sum_{[i,j]}sgn(i,j)\tilde R^{i_1i_2}_{j_1j_2}\dots\tilde R^{i_{N-1}i_{N}}_{j_{N-1}j_{N}}\right)|_{\pi(o)}dvol_M(o),
\end{eqnarray*}
where, as previously, $\hat R^{\alp\beta}_{\gamma\delta}|_o$ denotes the curvature of the fiber $\pi^{-1}\pi(o)$ at the point $o$, and $\tilde R^{ij}_{kl}|_{\pi(o)}$ denotes the curvature of the base $B$ at the point $\pi(o)$.

For $p^*\in B$ we denote by $Z_{p^*}:=\pi^{-1}(p^*)$, and by $dvol_{Z_{p^*}}$ the volume form on $Z_{p^*}$ with respect to the restriction of the original Riemannian metric $g$.
Integrating along the fibers we have
\begin{eqnarray*}
\lim_{\eps\to+0} V_{n-e}(M(\eps))=\\
(2\pi)^{-\frac{e-N}{2}}\int_B dvol_B(p^*)\left(\sum_{[\alp,\beta]}sgn(\alp,\beta) \hat R^{\alp_1\alp_2}_{\beta_1\beta_2}\dots  \hat R^{\alp_{e-N-1}\alp_{e-N}}_{\beta_{e-N-1}\beta_{e-N}}\right)|_{p^*} \cdot\\
(2\pi)^{-\frac{N}{2}}\int_{p\in Z_{p^*}}dvol_{Z_{p^*}}(p)
\left( \sum_{[i,j]}sgn(i,j)\tilde R^{i_1i_2}_{j_1j_2}\dots\tilde R^{i_{N-1}i_{N}}_{j_{N-1}j_{N}}\right)\big|_{p}
\end{eqnarray*}
The inner integral is equal to $\chi(Z_{p^*})$ by the generalized Gauss-Bonnet formula. Hence
\begin{eqnarray*}
\lim_{\eps\to+0} V_{n-e}(M(\eps))=\\
\chi(Z)\cdot (2\pi)^{-\frac{e-N}{2}}\int_B dvol_B(p^*)\left(\sum_{[\alp,\beta]}sgn(\alp,\beta) \hat R^{\alp_1\alp_2}_{\beta_1\beta_2}\dots  \hat R^{\alp_{e-N-1}\alp_{e-N}}_{\beta_{e-N-1}\beta_{e-N}}\right)\big|_{p^*}=\\
\chi(Z) V_{n-e}(B),
\end{eqnarray*}
where the last equality is by the definition of $V_{n-e}$. \qed


\section{When $M(\eps)$ has uniform lower bound on sectional curvature?}\label{S:proof2}
\begin{proposition}\label{P:sect-curvature}
The sectional curvature of $M(\eps)$ stays uniformly bounded below as $\eps\to +0$ if and only if the sectional curvature of all fibers of $\pi$ is non-negative.
\end{proposition}
{\bf Proof.} Below we use the convention of lifting and lowering indices using the matrix $(g_{pq})$ of size $\dim M$. For $p\ne q$ let us denote by $K_{pq}(\eps)$ the Gauss curvature of $g(\eps)$ of the $(x^px^q)$-plane.
By the known formula
$$K_{pq}(\eps)=\frac{2R_{pqpq}(\eps)}{g_{pp}(\eps)g_{qq}(\eps)-g_{pq}^2(\eps)}.$$

\underline{Case 1.} Assume $p=\alp, q=\beta$. By Lemma \ref{L:new-metric} 
\begin{eqnarray*}
K_{\alp \beta}(\eps)=\frac{2R_{\alp\beta\alp\beta}(\eps)}{g_{\alp\alp}g_{\beta\beta}-g_{\alp\beta}^2}\overset{(\ref{E:curv2})}{=}\frac{O(1)}{g_{\alp\alp}g_{\beta\beta}-g_{\alp\beta}^2}=O(1).
\end{eqnarray*}

\underline{Case 2.} Assume $p=\alp,q=i$. Similarly 
by Lemma \ref{L:new-metric}
\begin{eqnarray*}
K_{\alp i}(\eps)=\frac{2R_{\alp i\alp i}(\eps)}{g_{\alp\alp}g_{ii}-g_{\alp i}^2}\overset{(\ref{E:curv1})}{=}\frac{O(\eps)}{\eps g_{\alp\alp}g_{ii}-\eps^2g_{\alp i}^2}=O(1).
\end{eqnarray*}

\underline{Case 3.} Assume $p=i,q=j$. Again 
\begin{eqnarray}\label{E:case3}
K_{ij}(\eps)=\frac{2R_{ijij}(\eps)}{\eps^2(g_{ii}g_{jj}-g_{ij}^2)}.
\end{eqnarray}


We have 
\begin{eqnarray}\label{E:sectional-curv1}
2R_{ijij}(\eps)\overset{(\ref{E:RC-formula})}{=}\\\label{E:sectional-curv2}
2\pt_i\pt_j g_{ij}(\eps)-\pt_i\pt_i g_{jj}(\eps)- \pt_j\pt_j g_{ii}(\eps)+2(\Gamma_{tij}(\eps)\Gamma^t_{ij}(\eps)-\Gamma_{tii}(\eps)\Gamma^t_{jj}(\eps)).
\end{eqnarray}


The $\Gamma \Gamma$ summand can be rewritten more explicitly
\begin{eqnarray*}
\Gamma_{tij}(\eps)\Gamma_{sij}(\eps)g^{st}(\eps)-\Gamma_{tii}(\eps)\Gamma_{sjj}(\eps)g^{st}(\eps).
\end{eqnarray*}
Using Lemma \ref{L:matrices} for the matrix $G=g(\eps)$ the last expression is equal
\begin{eqnarray}\label{E:Christof-eps}
\Gamma_{kij}(\eps)\Gamma_{lij}(\eps)(\tilde g^{kl}(\eps)+O(1))-\Gamma_{kii}(\eps)\Gamma_{ljj}(\eps)(\tilde g^{kl}(\eps)+O(1)).
\end{eqnarray}
 We may choose furthermore the coordinates near $o^*$ to be normal.
Then by Lemma \ref{L:estimate-Christoffel} the expression (\ref{E:Christof-eps}) is equal to
\begin{eqnarray*}
\eps(\Gamma_{kij}\Gamma_{lij}\tilde g^{kl}-\Gamma_{kii}\Gamma_{ljj}\tilde g^{kl})+O(\eps^2).
\end{eqnarray*}
Substituting this back to (\ref{E:sectional-curv2}) instead of the $\Gamma\Gamma$ term we get
\begin{eqnarray}
2R_{ijij}(\eps)=\\
\eps(2\pt_i\pt_j g_{ij}-\pt_i\pt_i g_{jj}- \pt_j\pt_j g_{ii})+
\eps(\Gamma_{kij}\Gamma_{lij}\tilde g^{kl}-\Gamma_{kii}\Gamma_{ljj}\tilde g^{kl})+O(\eps^2).
\end{eqnarray}
Substituting this into (\ref{E:case3}) we get 
\begin{eqnarray*}
K_{ij}(\eps)=\frac{\eps(2\pt_i\pt_j g_{ij}-\pt_i\pt_i g_{jj}- \pt_j\pt_j g_{ii}+
\Gamma_{kij}\Gamma_{lij}\tilde g^{kl}-\Gamma_{kii}\Gamma_{ljj}\tilde g^{kl})+O(\eps^2)}{\eps^2(g_{ii}g_{jj}-g_{ij}^2)}=\\
\eps^{-1}\frac{2\tilde R_{ijij}}{g_{ii}g_{jj}-g_{ij}^2}+O(1)=\eps^{-1}\tilde K_{ij}+O(1),
\end{eqnarray*}
where the last summand $O(1)$ is uniform on $M$, and $\tilde K_{ij}$ denotes the sectional curvature of the fiber at $o$. We see that when $\eps\to +0$ the sectional curvature $K_{ij}(\eps)$ is uniformly on $M$ bounded below
if and only if $\tilde K_{ij}\geq 0$. \qed


\begin{thebibliography}{99}

\bibitem{alesker-mflds1}
Alesker, Semyon; Theory of valuations on manifolds. I. Linear spaces. Israel J. Math. 156 (2006), 311–339. 

\bibitem{alesker-mflds2}
 Alesker, Semyon; Theory of valuations on manifolds. II. Adv. Math. 207 (2006), no. 1, 420–454. 

\bibitem{alesker-mflds3}
Alesker, Semyon; Fu, Joseph H. G.; Theory of valuations on manifolds. III. Multiplicative structure in the general case. Trans. Amer. Math. Soc. 360 (2008), no. 4, 1951–1981.

\bibitem{alesker-mflds4}
Alesker, Semyon; Theory of valuations on manifolds. IV. New properties of the multiplicative structure. Geometric aspects of functional analysis, 1–44, Lecture Notes in Math., 1910, Springer, Berlin, 2007. 


\bibitem{alesker-kent}
Alesker, Semyon; Introduction to the theory of valuations. CBMS Regional Conference Series in Mathematics, 126. 
Published for the Conference Board of the Mathematical Sciences, Washington, DC; by the American Mathematical Society, Providence, RI, 2018.

\bibitem{alesker-conjectures}
Alesker, Semyon; Some conjectures on intrinsic volumes of Riemannian manifolds and Alexandrov spaces. Arnold Math. J. 4 (2018), no. 1, 1–17.

\bibitem{allendoerfer-weil}
Allendoerfer, Carl B.; Weil, Andr\'e;
The Gauss-Bonnet theorem for Riemannian polyhedra.
Trans. Amer. Math. Soc. 53 (1943), 101–129.


\bibitem{bernig-fu-unitary}
Bernig, Andreas; Fu, Joseph H. G.; Hermitian integral geometry. Ann. of Math. (2) 173 (2011), no. 2, 907–945.

\bibitem{bernig-fu-solanes-unitary}
Bernig, Andreas; Fu, Joseph H. G.; Solanes, Gil; Integral geometry of complex space forms. Geom. Funct. Anal. 24 (2014), no. 2, 403–492.

\bibitem{bernig-faifman}
Bernig, Andreas; Faifman, Dmitry; Valuation theory of indefinite orthogonal groups. J. Funct. Anal. 273 (2017), no. 6, 2167–2247.

\bibitem{bernig-faifman-solanes}
Bernig, Andreas; Faifman, Dmitry; Solanes, Gil; Curvature Measures of Pseudo-Riemannian Manifolds. arXiv:1910.09635

\bibitem{bernig-faifman-solanes2}
Bernig, Andreas; Faifman, Dmitry; Solanes, Gil;
Uniqueness of curvature measures in pseudo-Riemannian geometry. arXiv:2009.02230

\bibitem{ch-mu-sch}
Cheeger, Jeff; M\"uller, Werner; Schrader, Robert; On the curvature of piecewise flat spaces. Comm. Math. Phys. 92 (1984), no. 3, 405–454.

\bibitem{chern-1944}
 Chern, Shiing-shen; A simple intrinsic proof of the Gauss-Bonnet formula for closed Riemannian manifolds. Ann. of Math. (2) 45 (1944), 747–752.

\bibitem{fu-wannerer}
 Fu, Joseph H. G.; Wannerer, Thomas; Riemannian curvature measures. Geom. Funct. Anal. 29 (2019), no. 2, 343–381. 

\bibitem{hadwiger-thm}
Hadwiger, Hugo; {\itshape Vorlesungen \"uber Inhalt, Oberfl\"ache und Isoperimetrie.} (German) Springer-Verlag, Berlin-G\"ottingen-Heidelberg 1957.

\bibitem{klain-rota}
Klain, Daniel A.; Rota, Gian-Carlo; Introduction to geometric probability. Lezioni Lincee. [Lincei Lectures] Cambridge University Press, Cambridge, 1997.

\bibitem{landau-lifschitz}
Landau, L. D.; Lifshitz, E. M.; Course of theoretical physics, Vol. 2. The classical theory of fields. Fourth edition. 
Translated from the Russian by Morton Hamermesh. Pergamon Press, Oxford-New York-Toronto, Ont., 1975. 

\bibitem{lebedeva-petrunin}
Lebedeva, Nina; Petrunin, Anton; Curvature tensors on Alexandrov spaces. In preparation (2017).

\bibitem{petrunin-personal}
Petrinin, Anton; Personal discussions.

\bibitem{schneider-book} 
Schneider, Rolf; {\itshape Convex bodies: the Brunn-Minkowski theory.} Second expanded edition. 
Encyclopedia of Mathematics and its Applications, 151. Cambridge University Press, Cambridge, 2014. 

\bibitem{weyl-39}
Weyl, Hermann; On the Volume of Tubes. Amer. J. Math. 61 (1939),
no. 2, 461-472.



\end{thebibliography}
\end{document}